\documentclass[12pt]{article}

\usepackage{latexsym}
\usepackage{amsmath}
\usepackage{amssymb}
\usepackage{amsfonts}

\begin{document}
\parskip 0.2 cm
    \renewcommand{\theequation}{\arabic{equation}}

\def \o{\over}
\def \s{\sigma}
\def \Y\times U{\Gamma}
\def \g{\gamma}
\def\gm{\gamma(dy,du)}
\def\xm{\xi(dy,du)}
\def \a{\alpha}
\def \l{\lambda}
\def \disp{\displaystyle}
\def \up{\uparrow}
\def \U{{\cal U}}
\def \V{{\cal V}}
\def \K{{\cal K}}
\def \P{{\cal P}}
\def \D{{\cal D}}
\def \M{{\cal M}}
\def \Mp{{\cal M_+}}
\def \ve{\varepsilon}
\def\reals{I\!\!R}
\def \O{\Omega}
\def \tO{\tilde\Omega}
\def \ph{\varphi}
\def \D{\Delta}
\def \Dd{{\cal D}}
\def \by{\bar y}
\def \bu{\bar u}
\def \pl{\partial}
\def \dn{\downarrow}
\def\hf{\hfill{$\Box$}}
\def \bp{\bar\psi}
\def\be{\bar \eta}
\def\d{\delta}
\def \ve{\varepsilon}
\def \t{\theta}
\def \td{\tilde d}
\def \hd{\hat d}
\def \tk{\tilde k}
\def \G{\Gamma}

\newcommand{\R}{\mathbb{R}}

\newtheorem{Theorem}{Theorem}[section]
\newtheorem{Proposition}[Theorem]{Proposition}
\newtheorem{Remark}[Theorem]{Remark}
\newtheorem{Lemma}[Theorem]{Lemma}
\newtheorem{Corollary}[Theorem]{Corollary}
\newtheorem{Definition}[Theorem]{Definition}
\newtheorem{Example}[Theorem]{Example}

\newtheorem{theorem}{Theorem}[section]
\newtheorem{proposition}[Theorem]{Proposition}
\newtheorem{remark}[Theorem]{Remark}
\newtheorem{lemma}[Theorem]{Lemma}
\newtheorem{corollary}[Theorem]{Corollary}	
\newtheorem{definition}[Theorem]{Definition}

\begin{center}\Large{{\bf Linear Programming Estimates for Ces\`aro and  Abel  Limits of Optimal Values in Optimal Control Problems}}\end{center}

{\bf V. Gaitsgory$^a$ and I. Shvartsman$^b$}\\
 $^a$ {\it\small{Department of Mathematics, Macquarie University, Macquarie Park, NSW 2113, Australia } }\\
 $^b$ {\it \small{Department of Mathematics and Computer Science,
Penn State Harrisburg, Middletown, PA 17057, USA}}

\bigskip


{\bf Abstract.} We consider infinite horizon optimal control problems with time averaging and time discounting criteria and give estimates for the Ces\`aro and  Abel  limits of  their optimal values  in the case when they depend on the initial conditions. We establish that these limits are bounded from above by the optimal value of a certain infinite dimensional (ID) linear programming (LP) problem and that they are bounded from below by the optimal value of the corresponding dual problem.  (These estimates imply, in particular, that the Ces\`aro and  Abel limits  exist and are equal to each other if there is no duality gap). In addition, we obtain IDLP-based optimality conditions for the long run average optimal control problem, and we illustrate these conditions by an example. 

\section{Inroduction}

In this paper, we consider infinite horizon optimal control problems with time averaging and time discounting criteria and give estimates for the Ces\`aro and  Abel  limits of  their optimal values  in the case when they depend on the initial conditions. We establish that these limits are bounded from above by the optimal value of a certain infinite dimensional (ID) linear programming (LP) problem and that they are bounded from below by the optimal value of the corresponding dual problem.  (These estimates imply, in particular, that the Ces\`aro and  Abel limits  exist and are equal to each other if there is no duality gap). In addition, we obtain IDLP-based optimality conditions for the long run average optimal control problem.

An LP approach to  optimal control problems allows to use the convex duality theory and  LP-based numerical techniques for analysis and construction of optimal
 solutions, and this approach
has been extensively studied in both deterministic and stochastic settings.
For example,   LP formulations of problems of optimal control of deterministic systems considered on finite time and infinite time intervals were studied in \cite{Goreac-Serea,Her-Her-Lasserre,Lass-Trelat,Rubio,Vinter} and  \cite{GQ,GQ-1,GPS-2017,GPS-2018,QS}, respectively.

 Despite a great deal of attention that the development of the LP approach to control problems has
attracted, the issue about the validity of the LP representation for the Ces\`aro and  Abel  limits of  the optimal values in case these limits are dependent
on initial conditions  has not been addressed until recently. In fact, to the best of our knowledge, first results in this direction for deterministic controlled ODEs were obtained in \cite{BG}. This paper was followed up by \cite{BGS}, where stronger results  were obtained under weaker assumptions for deterministic controlled systems evolving in discrete time. (Note that  related results were obtained in \cite{HK-1}, \cite{HK-2} for Markov decision processes (MDP) with finite state/action spaces and in \cite{Gonzalez-Hernandes} for MDP with infinite uncountable state/action spaces).
In this paper, we use ideas from \cite{BGS} to strengthen and further develop some of the results of \cite{BG}.

Consider the  control  system
\begin{equation}\label{e-CSO}
y'(t)=f(y(t),u(t)), \ \ \ \ \ u(t)\in U, \ \ \ \ \ t\in [0,\infty) \ ,
\end{equation}
where $U$ is a compact metric space,
  $\ f(\cdot,\cdot):  \mathbb{R}^m\times U \to \reals^m$ is continuous
and satisfies Lipschitz condition in $y$ uniformly in $u \in U$. The controls
 $u(\cdot) $ are  measurable functions $u(\cdot): [0,\infty)\to U $, and the set of all such controls is denoted by $\U$. For a given  $u(\cdot) \in \U $ and initial condition $y(0)=y_0 $, the corresponding solution of (\ref{e-CSO}) is denoted by $y(t,y_0,u)$.

Let $Y\subset \reals^m$ be a non-empty compact set, equal to the closure of its interior.
We denote by $\U_T(y_0)$ and $\ \U(y_0)$  the sets of controls such that
 \begin{equation}\label{e-CSO-1}
y(t,y_0,u)\in Y \ \
\end{equation}
for any $\  t\in [0,T]$, (respectively, for any $\ t\in [0,\infty)$).
Inclusion (\ref{e-CSO-1})
can be interpreted as a state constraint. Everywhere  in what follows, it will be assumed the set $\U(y_0)$ is not empty for any $y_0\in Y$, that is, there exists at least one admissible control for any initial condition. (Systems that satisfy this property are called {\em viable} on $Y$.)

On the trajectories of system (\ref{e-CSO}), consider the  optimal control problem
\begin{equation}\label{Cesaro}
\frac{1}{T} \inf_{u(\cdot)\in \U_T(y_0)}\int_0^T k(y(t,y_0,u),u(t))dt=:V_T(y_0)
\end{equation}
and the optimal control problem
 \begin{equation}\label{Abel}
 \lambda \inf_{u(\cdot)\in \U(y_0)}\int_0^{\infty}e^{-\lambda t} k(y(t,y_0,u),u(t))dt=:h_{\lambda}(y_0),
\end{equation}
where $T>0 $, $\lambda>0 $ are positive parameters and
 $k(y,u):  \reals^m\times U \to \reals $ is a continuous function. Our focus will be on establishing estimates  for the
 limits $\disp \lim_{T\rightarrow\infty}V_T(y_0)$ and  $\disp \lim_{\l\dn 0}h_{\l}(y_0)$ (called
 Ces\`aro and  Abel  limits, respectively). Matters related to the existence and the equality of  Ces\`aro and  Abel limits of the optimal values have been addressed by many authors (see, e.g.,
\cite{Arisawa-3,Bardi,BQR-2015,GQ,GruneSIAM98,GruneJDE98,Khlopin,Sorin92,OV-2012,QR-2012}). A special feature and the novelty of our consideration is that we are making use of
occupational measure reformulations of  problems (\ref{Cesaro}) and (\ref{Abel}) and utilize the LP duality theory.

 The paper is organized as follows. Section 2 is devoted to preliminaries. In particular, in this section, we reformulate optimal control problems (\ref{Cesaro}) and (\ref{Abel})
 in terms of occupational measures, and  we introduce the IDLP problem and its dual that are instrumental for our consideration.
In Section 3 we establish that $\disp \limsup_{T\rightarrow\infty}V_T(y_0)$ and  $\disp \limsup_{\l\dn 0}h_{\l}(y_0)$ are less or equal than the optimal value of this IDLP problem and that
$\disp \liminf_{T\rightarrow\infty}V_T(y_0)$ and  $\disp \liminf_{\l\dn 0}h_{\l}(y_0)$ are greater or equal than the optimal value of its dual. (A straightforward corollary of these statements is the fact that the  Ces\`aro and  Abel  limits exist  and are equal if there is no duality gap).
In Section 4 we discuss IDLP-based optimality conditions for the long run average optimal control problem. In Section 5, we consider an example illustrating results obtained in Sections 3 and 4.
Some of the results of this paper were announced in \cite{CDC2019}, but, in contrast to \cite{CDC2019} (and to \cite{BG}), we do not assume in this paper that the set $Y$ if forward invariant; we assume only that the system is viable on $Y$.

\section{Occupational Measures and IDLP Problems}

Denote by
 $\ \mathcal{P}(Y\times U)$,  $\ \mathcal{M}_+(Y\times U)$ and $\ \mathcal{M}(Y\times U)$ the spaces of probability measures, non-negative measures and all finite measures, respectively, defined on the Borel subsets of $Y\times U$. The convergence in these  spaces is understood in the weak$^*$ sense,  that is, $\gamma^k \in \mathcal{M}(Y\times U), k =1,2,... ,$ converges to $\gamma
\in \mathcal{M}(Y\times U)$ if and only if
\begin{equation}\label{weak}
 \lim_{k\rightarrow \infty}\int_{Y\times U} \phi(y,u) \gamma^k (dy,du) \ = \
 \int_{Y\times U} \phi(y,u) \gamma (dy,du)
 \end{equation}
for any continuous $\phi(y,u): Y\times U \rightarrow \reals$.

Let $u(\cdot) \in \U_T(y_0)$ and  $y(t) =
y(t,y_0,u(\cdot)), \ t\in [0,T] $. A probability measure
 $\gamma_{u(\cdot),T} \in {\cal P} (Y \times U)$ is called the
 {\it occupational measure} generated by the pair $(y(\cdot),u(\cdot) )$ on the interval $[0,T]$ if, for
  any Borel set $Q \subset Y \times
 U$,
 $$
  \gamma _{u(\cdot),T} (Q) = \frac{1}{T}\int _0 ^T
 1_Q (y(t),u(t)) dt ,
 $$
 where $1_Q (\cdot)$ is
the indicator function of $Q$. It can be shown that  this definition is equivalent to
 the equality
\begin{equation}\label{e:occup-meas-def-eq-S}
\int_{Y\times U} q(y,u)\gamma_{u(\cdot),T} (dy,du) = \frac{1}{T} \int _0 ^
T q (y(t),u(t)) dt
 \end{equation}
 for any
$q(\cdot)\in C(Y\times U)$ (see, e.g., Section 1 in \cite{GPS-2018}).

Let $u(\cdot) \in {\cal U}(y_0)$  and  $y(t)=y(t,y_0,u(\cdot)),
\ t\in [0,\infty)  $. A probability measure
 $\gamma^{\lambda}_{u(\cdot)} \in {\cal P} (Y \times U)$  is called the {\it discounted occupational measure} generated
 by the pair $(y(\cdot),u(\cdot) )$ if for any Borel set $Q \subset Y \times
 U$,
$$
\gamma ^{\lambda}_{u(\cdot)} (Q) = \lambda \int _0 ^ \infty
e^{-\lambda t} 1_Q (y(t),u(t)) dt.
$$
The latter definition
is  equivalent to  the equality
\begin{equation}\label{e:occup-meas-def-eq}
\int_{Y\times U} q(y,u)\gamma ^{\lambda}_{u(\cdot)} (dy,du) = \lambda \int _0 ^ \infty
e^{-\lambda t} q (y(t),u(t)) dt
\end{equation}
for any
$q(\cdot)\in C(Y\times U)$.

Let $\Gamma_T(y_0)$ and $\Theta^{\lambda}(y_0) $ denote the sets of all occupational and discounted occupational measures, respectively, generated by
admissible processes:
$$
\Gamma_T(y_0):= \bigcup_{u(\cdot)\in\U_T(y_0)}\{\gamma _{u(\cdot),T}\}, \;
\Theta^{\lambda}(y_0):= \bigcup_{u(\cdot)\in\U(y_0)}\{\gamma^{\lambda}_{u(\cdot)}\}
$$
and
\begin{equation}\label{PZ1}
\Gamma_T:=\bigcup_{y_0\in Y}\Gamma_T(y_0),\quad  \Theta^{\lambda}:=\bigcup_{y_0\in Y}\Theta^{\lambda}(y_0).
\end{equation}

Due to (\ref{e:occup-meas-def-eq-S}) and (\ref{e:occup-meas-def-eq}), problems (\ref{Cesaro}) and (\ref{Abel}) can be equivalenlty reformulated as
\begin{equation}\label{e:occup-meas-def-eq-2}
\inf_{\gamma\in \Gamma_T(y_0) }\int_{Y\times U}k(y,u)\gamma(dy,du) =V_T(y_0),
\end{equation}
and
\begin{equation}\label{e:occup-meas-def-eq-3}
\inf_{\gamma\in \Theta^{\lambda}(y_0) }\int_{Y\times U}k(y,u)\gamma(dy,du) =h_{\l}(y_0).
\end{equation}
To describe convergence properties of occupational measures, we introduce the following metric on $\P({Y\times U})$:
$$
\rho(\g',\g''):=\sum_{j=1}^{\infty} {1\o 2^j}\left|\int_{Y\times U} q_j(y,u)\g'(dy,du)-\int_{Y\times U} q_j(y,u)\g''(dy,du)\right|
$$
for $\g',\g''\in \P({Y\times U})$, where $q_j(\cdot),\,j=1,2,\dots,$ is a sequence of Lipschitz continuous functions dense in the unit ball of the space of continuous functions $C({Y\times U})$ from ${Y\times U}$ to $\reals$.
This metric is consistent with the weak$^*$ convergence topology on $\P({Y\times U})$, that is,
a sequence $\g^k\in \P({Y\times U})$ converges to $\g\in \P({Y\times U})$ in this metric if and only if \eqref{weak} holds.

Using the metric $\rho$, we can define the ``distance" $\rho(\g,\Gamma)$ between $\g\in \P({Y\times U})$ and $\Gamma\subset \P({Y\times U})$
and the Hausdorff metric $\rho_H(\Gamma_1,\Gamma_2)$ between $\Gamma_1\subset \P({Y\times U})$ and $\Gamma_2\subset \P({Y\times U})$ as follows:
$$
\rho(\g,\Gamma):=\inf_{\g'\in \Gamma}\rho(\g,\g'),\quad
\rho_H(\Gamma_1,\Gamma_2):=\max\{\sup_{\g\in \Gamma_1}\rho(\g,\Gamma),\sup_{\g\in \Gamma_2}\rho(\g,\Gamma_2)\}.
$$

Define the set $W$ by
\begin{equation*}
\begin{aligned}
W:=\{\g\in \P(Y\times U)|\, \int_{Y\times U}& \nabla \phi(y)f(y,u)\,d\gamma(dy,du)=0
\hbox{ for all }\phi\in C^1\}.
\end{aligned}
\end{equation*}
Here and below gradients are understood as {\em row} vectors, and $C^1$ is the space of continuously differentiable functions on $\reals^n$.

Consider the IDLP problem
\begin{equation}\label{limits-non-ergodic}
\inf_{(\gamma, \xi)\in \Omega(y_0)}\int_{Y\times U}k(y,u)\gamma(dy,du)=: k^*(y_0),
\end{equation}
where
$$
 \Omega(y_0):= \{(\gamma, \xi)\in \mathcal{P}(Y\times U)\times \mathcal{M}_{+}(Y\times U) |\,\gamma\in W,
 $$
\begin{equation}\label{non-ergodic-Omega}
\begin{aligned}
&\int_{Y\times U}(\phi(y_0)-\phi(y))\gamma(dy,du) +
 \int_{Y\times U}\nabla \phi(y)f(u,y)\xi(dy,du)    =0 \ \ \forall \phi(\cdot)\in C^1 \}.
\end{aligned}
\end{equation}
(This is an LP problem since both the objective function and the constraints are linear in the ``decision variables" $\g$ and $\xi$.)

Consider also the IDLP problem
\begin{equation}\label{limits-non-ergodic-dual}
 \sup_{(\mu , \psi(\cdot), \eta(\cdot) )\in \mathcal{D}}\mu :=d^*(y_0),
\end{equation}
where $\mathcal{D}$ is the set of triplets $(\mu , \psi(\cdot), \eta(\cdot) )\in \reals^1\times C^1\times C^1$
that for all $(y,u)\in Y\times U$ satisfy the inequalities
\begin{equation}\label{limits-non-ergodic-dual-1}
k(y,u)+ (\psi (y_0)- \psi (y)) + \nabla \eta (y) f(y,u)-\mu \geq 0,
\end{equation}
\begin{equation}\label{limits-non-ergodic-dual-2}
 \nabla \psi (y) f(y,u)\geq 0.
\end{equation}
The optimal value of problem (\ref{limits-non-ergodic-dual})  can be equivalently represented as
\begin{equation}\label{CC9}
d^*(y_0)=\sup_{(\psi,\eta)}\min_{(y,u)}\{k(y,u)+(\psi(y_0)-\psi(y))+\nabla\eta(y)f(y,u)\}.
\end{equation}

Problem (\ref{limits-non-ergodic-dual})-(\ref{limits-non-ergodic-dual-2}) can be shown to be dual to (\ref{limits-non-ergodic}) (see Section 5 of \cite{BG}) and it is proved in \cite{BG}, Theorem 3.1,  that
\begin{equation}\label{limits-non-ergodic-dual-4}
 k^*(y_0) \geq d^*(y_0).
\end{equation}

\bigskip

In this paper we will be using Clarke's generalized gradient of a Lipschitz function $\phi:\,\reals^n\to \reals$, which  can be represented as a convex hull of the limits of the gradients, that is, (see, e.g.,\cite{Bardi}, p.63)
\begin{equation}\label{KH1}
\partial \phi(y)={\rm co}\,\{p|\,p=\lim_{k\to \infty}\nabla \phi(y_k)\hbox{ for some }y_k\to y\}.
\end{equation}

The following approximation property will also be used.
\begin{proposition}\label{subd} (\cite{CZAR}, Theorem 2.2)
Let $E$ be an open subset of $\reals^n$ and let $\phi:\,E\to \reals$ be a locally Lipschitz function. Then for any $\ve>0$ there exists a function $\phi_{\ve}:\,E\to \reals$ of class $C^{\infty}$ such that for any $y\in E$
$$
|\phi_{\ve}(y)-\phi(y)|\le \ve
$$
and
$$
\nabla \phi_{\ve}(y)\in \bigcup_{y'\in (y+\ve B)\cap E}\pl \phi(y')+\ve B,
$$
where $B$ is an open unit ball.
\end{proposition}

\section{Estimates for the  Ces\`aro and  Abel  Limits of the Optimal Values}

The estimates from below are established by the following proposition.

\begin{proposition}\label{Prop-former-2-3}
The following inequalities hold:
\begin{equation}\label{e-main-1}
\liminf_{T\rightarrow\infty}V_T(y_0)\geq d^*(y_0),
\end{equation}
and
\begin{equation}\label{e-main-2}
\liminf_{\lambda\downarrow 0}h_{\lambda}(y_0)\geq d^*(y_0).
\end{equation}
\end{proposition}
These inequalities were proved in Proposition 3.1 in \cite{BG} under the assumption of invariance of $Y$, which is not used in the proof and can be replaced with viability. The proof below is much simpler than in \cite{BG}.

{\bf Proof of Proposition \ref{Prop-former-2-3}.}
Assume that \eqref{e-main-1} is not true, that is, $\disp \liminf_{T\to \infty} V_T(y_0)< d^*(y_0)$. Then there exist  $\ \beta >0$ and a pair of $C^1$ functions $\ (\psi(\cdot),\eta(\cdot))$ with $\ \psi(\cdot) $ satisfying \eqref{limits-non-ergodic-dual-2}, such that
$$
k(y,u) + (\psi(y_0)-\psi(y)) +\nabla\eta(y)f(y,u)\geq \liminf_{T'\to \infty} V_{T'}(y_0) +\beta.
$$
Taking into account that $\psi(y(t))$ is non-decreasing along an arbitrary admissible process $(y(\cdot),u(\cdot)) $ due to \eqref{limits-non-ergodic-dual-2}, we have from the last inequality that
$$
k(y(t),u(t)) + \nabla\eta(y(t))f(y(t),u(t))\geq \liminf_{T'\to \infty} V_{T'}(y_0) +\beta, \ \ \ t\ge 0.
$$
Therefore,
$$
 \frac{1}{T}\int_{0}^{T}k(y(t),u(t))\,dt + \frac{1}{T}(\eta(y(T))- \eta(y_0))\geq \liminf_{T\to \infty} V_{T'}(y_0) +\beta,
$$
which implies that
$$
 V_T(y_0) + \frac{1}{T}(\max_{y\in Y}\eta(y)- \eta(y_0))\geq \liminf_{T'\to \infty} V_{T'}(y_0) +\beta.
$$
By taking $\liminf_{T\rightarrow\infty}$ on the left side, we obtain a contradiction. Thus, \eqref{e-main-1} is proved.

Similarly, to prove \eqref{e-main-2} assume that it is not true and that there exist  $\ \beta >0$ and a pair of $C^1$ functions $\ (\psi(\cdot),\eta(\cdot))$ with $\ \psi(\cdot) $ satisfying \eqref{limits-non-ergodic-dual-2}, such that
$$
k(y,u) + (\psi(y_0)-\psi(y)) +\nabla\eta(y)f(y,u)\geq \liminf_{\l'\dn 0} h_{\l'}(y_0) +\beta.
$$
This implies that for an arbitrary admissible process $(y(\cdot),u(\cdot))$ and $\l>0$
\begin{equation}\label{VM1}
\l \int_{0}^{\infty}e^{-\l t}k(y(t),u(t))\,dt + \l\int_{0}^{\infty}e^{-\l t}\nabla\eta(y(t))f(y(t),u(t))\,dt\geq  \liminf_{\l'\dn 0} h_{\l'}(y_0) +\beta.
\end{equation}
Taking into account that
\begin{equation*}
\begin{aligned}
& \l\int_{0}^{\infty}e^{-\l t}\nabla\eta(y(t))f(y(t),u(t))\,dt= \l\int_{0}^{\infty}e^{-\l t}{d\o dt}\eta(y(t))\,dt\\
&=-\l\eta(y_0)+\l^2\int_{0}^{\infty}e^{-\l t}\eta(y(t))\,dt\to 0
\hbox{ as }\l\dn 0,
\end{aligned}
\end{equation*}
we conclude from \eqref{VM1} that
$$
\liminf_{\lambda\downarrow 0}h_{\lambda}(y_0)\ge \liminf_{\lambda'\downarrow 0}h_{\lambda}(y_0)+\beta,
$$
which is a contradiction. The proposition is proved. \hf

\bigskip

Let us now establish the estimate from above. To this end, take $\d>0$ and denote:
$$
Y^{\d}:=Y+\d \bar B,
$$
where $\bar B$ is the closed unit ball.
For $y_0\in Y^{\d}$, along with problems \eqref{Cesaro} and \eqref{Abel}, consider the problems
\begin{equation}\label{Cesaro-1}
\frac{1}{T} \inf_{u(\cdot)\in \U^{\d}_T(y_0)}\int_0^T k(y(t,y_0,u),u(t))dt=:V^{\d}_T(y_0),
\end{equation}
\begin{equation}\label{Abel-1}
\lambda \inf_{u(\cdot)\in \U^{\d}(y_0)}\int_0^{\infty}e^{-\lambda t} k(y(t,y_0,u),u(t))dt=:h^{\d}_{\lambda}(y_0),
\end{equation}
where $\U^{\d}_T(y_0)\subset \U$ and $\U^{\d}(y_0)\subset \U$ are the sets of controls such that, for any $u(\cdot)\in \U^{\d}_T(y_0)$ (respectively, $u(\cdot)\in \U^{\d}(y_0)$), the corresponding trajectories satisfies the inclusions $y(t,y_0,u(\cdot))\in Y^{\d} \hbox{ for all }t\in[0,T]$ (respectively, $y(t,y_0,u(\cdot))\in Y^{\d} \hbox{ for all }t\ge 0$).

We will be referring to the following assumptions for all $\d\in(0,\d_0)$ for some $\d_0>0$:

(A1) The set $Y^{\d}$ is viable, that is, $\U^{\d}(y_0)\neq \emptyset$ for all $y_0\in Y^{\d}$.

(A2) The function $y_0\mapsto V^{\d}_T(y_0)$ is Lipschitz on $Y^{\d}$ for all $T>0$.

(A3) The function $y_0\mapsto h^{\d}_{\l}(y_0)$ is Lipschitz on $Y^{\d}$ for all $\l>0$.

(A4) For any $y_0\in Y$
\begin{equation}\label{MC2}
\lim_{T\to \infty,\d\to 0}|V_T(y_0)-V_T^{\d}(y_0)|=0.
\end{equation}

(A5) There holds the relation
\begin{equation}\label{GH2}
\lim_{T\to \infty}\rho_H({\rm co}\, \Gamma_T,W)=0.
\end{equation}

(A6) For any $y_0\in Y$
\begin{equation}\label{MC3}
\lim_{\l\dn 0,\d\to 0}|h_{\l}(y_0)-h_{\l}^{\d}(y_0)|=0.
\end{equation}

(A7) There holds the relation
\begin{equation}\label{MC4}
\lim_{T\to \infty}\rho_H({\rm co}\, \Theta^{\l},W)=0.
\end{equation}
(Recall that $\Gamma_T$ and $\Theta^{\l}$ in (A5) and (A7) are defined in \eqref{PZ1}.)

\medskip

\begin{remark}
{\rm A sufficient condition for (A5) to be valid is that, for any Lipschitz continuous function $q(y,u):\, \reals^n\times \reals^m\to \reals$,
\begin{equation}\label{MC1}
{1\o T}\left|\inf_{u(\cdot)\in \U_T(y_0), y_0\in Y }\int_0^T q(y(t,y_0,u),u(t))\,dt-\inf_{u(\cdot)\in \U^{\d}_T(y_0),  y_0\in Y }\int_0^T q(y(t,y_0,u),u(t))\,dt\right|\to 0
\end{equation}
as $\d\to 0$ and $T\to \infty$; see Theorem 2.1 in \cite{G04}. A sufficient condition for (A7) to be valid is that, for any Lipschitz continuous function $q(y,u):\, \reals^n\times \reals^m\to \reals$ and any $y_0\in Y $,
\begin{equation}\label{MC1-1}
\lambda\left|\inf_{u(\cdot)\in \U(y_0) }\int_0^{\infty} e^{-\lambda t} q(y(t,y_0,u),u(t))dt-\inf_{u(\cdot)\in \U^{\d}(y_0) }\int_0^{\infty} e^{-\lambda t} q(y(t,y_0,u),u(t))dt\right|\to 0
\end{equation}
as $\d\to 0$ and $\lambda\to 0$, with the convergence being uniform with respect to $y_0\in Y$. The fact that this condition  implies (A7) can be established similarly
 to Proposition 6.1 in \cite{GQ}. Note that (\ref{MC1}) and (\ref{MC1-1}) are stronger versions of the assumptions (A4) and (A6) (respectively). Note also that
 both (\ref{MC1}) and (\ref{MC1-1})
are satisfied if $Y $ is invariant with respect to the solutions of system (\ref{e-CSO}), that is, if all trajectories of (\ref{e-CSO}) with initial condition in $Y$ do not leave $Y$ (since  $\U^{\d}(y_0)=\U(y_0)=\U$ for any $y_0\in Y$ in this case).

It is worth noting that (A5) and (A7) are satisfied automatically in the relaxed control setting, which is a consequence of Theorems 2.2, 2.5 and Lemma 2.4 in \cite{GQ-1}.
}
\end{remark}

\begin{Theorem}\label{Th-upper-bound}
{\rm (a)} Under assumptions {\rm (A1), (A2), (A4),} and {\rm(A5)},
\begin{equation}\label{ub-1}
\limsup_{T\rightarrow\infty} V_{T}(y_0)\leq k^*(y_0) \ \ \forall \ y_0\in Y.
\end{equation}
{\rm (b)} Under assumptions {\rm (A1), (A3), (A6),} and {\rm(A7)},
\begin{equation}\label{ub-2}
\limsup_{\l\dn 0} h_{\l}(y_0)\leq k^*(y_0) \ \ \forall \ y_0\in Y.
\end{equation}
\end{Theorem}

\bigskip

To prove this theorem we need a few auxiliary results. First of these is the following lemma. (This lemma constitutes a part of the proof of Theorem 3.8 in \cite{BQR-2015}; its proof  is omitted in \cite{BQR-2015}, and we prove it here since it is important for our consideration.)

\begin{lemma} Under assumption (A5),
\begin{equation}\label{e-before-lim-1}
\int_{Y\times U} V_T(y)\,\g(dy,du)\le   \int_{Y\times U} k(y,u)\,\g(dy,du)\ \ \ \forall \ \gamma\in W
\end{equation}
for any $T>0$.
\end{lemma}
{\bf Proof.}
Let us show first that, for any $T,T'>0$
\begin{equation}\label{CC20}
\int_{Y\times U} V_T(y)\,\g'(dy,du)\le   \int_{Y\times U} k(y,u)\,\g'(dy,du)+{MT\o T'}\quad \ \forall \ \g'\in  \G_{T'}(y_0), \ \forall \ y_0\in Y,
\end{equation}
where $M$ is such that $|f(y,u)|\le M$ for all $(y,u)\in Y\times U$.

Take $y_0\in Y$, $\g'\in \G_{T'}(y_0)$, and let  $u(\cdot)\in \U_{T'}(y_0)$ be a control that generates $\g'$ on $[0,T']$. Extend $u$ from the interval $[0,T]$ to the interval $[0,T'+T]$ so that $u\in \U_{T'+T}(y_0)$. (Such extension is possible due to viability of $Y$.) Let $y(\cdot)$ be the corresponding trajectory. Taking into account that $\disp V_T(y(s))\le  {1\o T}\int_{0}^{T} k(y(r+s),u(r+s)))\,dr$ for all $s\in [0,T']$, we obtain
\begin{equation*}
\begin{aligned}
&\int_{Y\times U} V_T(y)\,\g'(dy,du)={1\o T'} \int_{0}^{T'} V_T(y(s))\,ds
\le {1\o T'} \int_{0}^{T'}{1\o T} \int_{0}^{T} k(y(r+s),u(r+s)))\,dr\,ds\\
&={1\o T} \int_{0}^{T}{1\o T'} \int_{0}^{T'} k(y(r+s),u(r+s))\,ds\,dr\\
&={1\o T} \int_{0}^{T}{1\o T'} \int_{r}^{T'+r} k(y(\s),u(\s))\,d\s\,dr
\le {1\o T} \int_{0}^{T}{1\o T'} \left(\int_{0}^{T'} k(y(\s),u(\s))+2Mr\right)\,d\s\,dr\\
&= {1\o T} \int_{0}^{T}{1\o T'} \int_{0}^{T'} k(y(\s),u(\s))\,d\s\,dr+ {1\o TT'} \int_{0}^{T}2Mr\,dr\\
&={1\o T} \int_{0}^{T}\int_{Y\times U} k(y,u)\,\g'(dy,du)\,dr+{MT\o T'}=\int_{Y\times U} k(y,u)\,\g'(dy,du)+{MT\o T'}\,.
\end{aligned}
\end{equation*}
Thus, inequality (\ref{CC20}) is established. From this inequality
it follows that
\begin{equation}\label{CC20-2}
\int_{Y\times U} V_T(y)\,\g'(dy,du)\le   \int_{Y\times U} k(y,u)\,\g'(dy,du)+{MT\o T'}\quad \ \ \ \forall\g'\in {\rm co}\ \G_{T'}\ .
\end{equation}
Due to (A5), for any $\gamma\in W $ there exist sequences $T'_l > 0, \ \g'_l\in {\rm co}\,\G_{T'_l} $, $\ l=1,2,...,$ such that $T'_l\rightarrow\infty $ and $\g'_l\rightarrow \g $. Passing to the limit along these sequences in (\ref{CC20-2})
and taking into account that
$$
\int_{Y\times U} V_T(y)\,\g(dy,du)=\lim_{\g'_l\rightarrow \g }\int_{Y\times U} V_T(y)\,\g'_l(dy,du)
$$
due to continuity of $V_T (\cdot)$, we arrive at  inequality (\ref{e-before-lim-1}). The lemma is proved. \hf

\bigskip

The following corollary follows from the fact that $V_T^{\d}(y)\le V_T(y)$ for all $y\in Y$.

\begin{corollary}
For any $T>0$ and any $\d\in (0,\d_0)$ we have
\begin{equation}\label{e-before-lim-10}
\int_{Y\times U} V_T^{\d}(y)\,\g(dy,du)\le   \int_{Y\times U} k(y,u)\,\g(dy,du)\ \ \ \forall \ \gamma\in W.
\end{equation}
\end{corollary}

Recall that, due to Radamacher's theorem, a Lipschitz function is almost everywhere differentiable.

\begin{Lemma}\label{L-DPP}
Let $V^{\d}_T(\cdot)$ be differentiable at $y_0\in {\rm int}\,Y^{\d}$ for some $T>0$ and $\d>0$. Then
\begin{equation}\label{ZZ1}
\nabla V^{\d}_T(y_0)f(y_0,u)\ge -{2M\o T} \quad\hbox{for all } u\in U,
\end{equation}
where $M$ is such that $|k(y,u)|\le M$ for all $(y,u)\in Y^{\d_0}\times U$.
\end{Lemma}

{\bf Proof.}
Let $u(\cdot)\in \U$ be a constant function, and $y(\cdot)=y(\cdot,y_0,u)$ be the corresponding trajectory.
From the dynamic programming principle it follows that  for any sufficiently small $\D t\in (0,T)$ such that $y(t)\in Y^{\d}$
for $t\in [0,\D t]$,
\begin{equation}\label{ZZ2}
\begin{aligned}
&TV^{\d}_T(y_0)\le
& \int_0^{\D t} k(y(t),u)\,dt+(T-\D t)V^{\d}_{T-\D t}(y(\D t)).
\end{aligned}
\end{equation}
Taking into account that
$$
\int_0^{\D t} k(y(t),u)\,dt=k(y_0,u)\D t+o(\D t)
$$
and
\begin{equation*}
\begin{aligned}
(T-\D t)V^{\d}_{T-\D t}(y(\D t))\le TV^{\d}_{T}(y(\D t))+M\D t,
\end{aligned}
\end{equation*}
we obtain
\begin{equation}\label{XX3}
TV^{\d}_T(y_0)\le k(y_0,u)\D t+TV^{\d}_T(y(\D t))+M\D t+o(\D t).
\end{equation}
Under the differentiability assumption,
$$
V^{\d}_T(y(\D t))=V^{\d}_{T}(y_0)+\nabla V^{\d}_{T}(y_0)f(y_0,u)\D t +o_T(\D t).
$$
Substituting this into \eqref{XX3} and passing to the limit as $\D t\to 0$ we obtain
\begin{equation*}
\nabla V^{\d}_T(y_0)f(y_0,u)\ge -{2M\o T}.
\end{equation*}
The lemma is proved.
\hf

\bigskip

{\bf Proof of Theorem \ref{Th-upper-bound}.} We will prove only the statement (a) of the theorem. The proof of the statement (b) is similar (see the Remark after the end of the proof). Due to \eqref{KH1} and Lemma \ref{L-DPP}, we have that
for any $\d\in (0,\d_0)$, $(y,u)\in  ({\rm int}\,Y^{\d})\times U$, and $T>0$
\begin{equation}\label{ZZ7}
pf(y,u)\ge -{2M\o T}\quad\hbox{for any }p\in \pl V^{\d}_T(y).
\end{equation}
(Here and below generalized gradients are {\em row} vectors.)
Due to Proposition \ref{subd} with $E={\rm int}\,Y^{\d}$, for any $\ve\in(0,\d)$ there exists $\psi_{\ve,\d,T}\in C^1({\rm int}\,Y^{\d})$ such that
\begin{equation}\label{ZZ6}
\max_{y\in Y} |\psi_{\ve,\d,T}(y)-V_T^{\d}(y)|\le \ve,
\end{equation}
and
\begin{equation}\label{ZZ8}
\nabla \psi_{\ve,\d,T}(y)\in \bigcup_{y'\in y+\ve B}\pl V_T^{\d}(y')+\ve B\quad\hbox{ for all }y\in Y,
\end{equation}
where $B$ is the open unit ball. From \eqref{ZZ7} and \eqref{ZZ8} we deduce that
there exists a constant $M$ such that
\begin{equation}\label{ZZ9}
\nabla \psi_{\ve,\d,T}(y)f(y,u)\ge -{2M\o T}-\ve M \quad\hbox{ for all }(y,u)\in Y\times U.
\end{equation}
Due to \eqref{ZZ6} we can rewrite \eqref{e-before-lim-10} as
\begin{equation*}
\int_{Y\times U}(k(y,u)+\ve-\psi_{\ve,\d,T}(y))\,\gm\ge 0\quad \hbox{for all }\g\in W,
\end{equation*}
which, in turn, is equivalent to
\begin{equation}\label{CC15-1}
\min_{\g\in W}\int_{Y\times U}(k(y,u)+\ve-\psi_{\ve,\d,T}(y))\,\gm\ge 0.
\end{equation}
The problem on the left hand side of (\ref{CC15-1}), i.e.,
\begin{equation}\label{CC8-1}
\min_{\g\in W}\int_{Y\times U}(k(y,u)+\ve-\psi_{\ve,\d,T}(y))\,\gm,
\end{equation}
is an IDLP problem. It is proved in \cite{GQ}, Section 3, that its dual is
\begin{equation}\label{BB12-1}
\sup_{\eta\in C^1}\min_{(y,u)\in Y\times U}\{k(y,u)+\ve-\psi_{\ve,\d,T}(y)+\nabla \eta(y)f(y,u)\}.
\end{equation}
The optimal values of \eqref{CC8-1} and \eqref{BB12-1} are equal (\cite{GQ}, Theorem 3.1). Therefore, \eqref{CC15-1} is equivalent to
\begin{equation}\label{BB15-1}
\sup_{\eta\in C^1}\min_{(y,u)\in Y\times U}\{k(y,u)+\ve-\psi_{\ve,\d,T}(y)+\nabla \eta(y)f(y,u)\}\ge 0.
\end{equation}
From \eqref{BB15-1} it follows that there exists a function $\eta_{\ve,\d,T}(\cdot)\in C^1$ such that for all $(y,u)\in Y\times U$
\begin{equation}\label{e-eps-feasib-1}
k(y,u)+\ve-\psi_{\ve,\d,T}(y)+\nabla\eta_{\ve,\d,T}f(y,u)\ge -\ve.
\end{equation}

For arbitrary $T>0$ and $\ve>0$ consider the following IDLP problem:
\begin{equation}\label{BB21-1}
\sup_{(\psi,\eta)\in Q(\ve,T)} \psi(y_0)=:  d^*(\ve,T,y_0),
\end{equation}
where
$Q(\ve,T)$ is the set of pairs $(\psi(\cdot),\eta(\cdot))\in C^1\times C^1$ that satisfy the inequalities
\begin{equation}\label{BB25-1}
\begin{aligned}
&k(y,u)-\psi(y)+\nabla \eta(y)f(y,u)\ge -2\ve,\\
&\nabla \psi(y)f(y,u)\ge -\frac{2M}{T}-M\ve\ \ \ \ \forall \ (y,u)\in Y\times U.
\end{aligned}
\end{equation}
From \eqref{ZZ9} and \eqref{e-eps-feasib-1} it follows that $\psi_{\ve,\d,T}$ and $\eta_{\ve,\d,T}$ satisfy \eqref{BB25-1}. ($\psi_{\ve,\d,T}$ in \eqref{ZZ7}-\eqref{ZZ6} is defined on int$\,Y^{\d}$, but it can be extended to $\reals^n$ without changing values on $Y$; therefore, without loss of generality, we can assume that $\psi_{\ve,\d,T}\in C^1.$) Therefore, taking into account \eqref{ZZ6},
\begin{equation}\label{ZZ10}
V_T^{\d}(y_0)-\ve\le \psi_{\ve,\d,T}(y_0)\le d^*(\ve,T,y_0).
\end{equation}
Similarly to \cite{BGS}, Lemma 6.1, it can be shown that
\begin{equation}\label{e-inequality-lemma-1}
 d^*(\ve,T, y_0)\leq k^*(\ve,T, y_0),
\end{equation}
where
\begin{equation*}
\begin{aligned}
k^*(\ve,T, y_0)=\inf_{(\g,\xi)\in \O(y_0)}\Big\{ \int_{Y\times U} k(y,u)\gm
+\left( \frac{2M}{T}+M\ve\right)\int_{Y\times U}\xi(dy,du)\Big\}.
\end{aligned}
\end{equation*}

It follows from \eqref{ZZ10} and \eqref{e-inequality-lemma-1} that
\begin{equation}\label{GH3}
V_T^{\d}(y_0)-\ve\le k^*(\ve,T,y_0).
\end{equation}
Next we establish that
\begin{equation}\label{BB1-T-2}
\lim_{\ve\downarrow 0,\,T\rightarrow\infty}k^*(\ve,T, y_0)= k^*(y_0).
\end{equation}
It is clear that $k^*(\ve,T, y_0) $ decreases as $T\to \infty$ and $\ve\dn 0$, and that $\ k^*(\ve,T, y_0) \geq k^*( y_0)   $ for any $T>0$ and $\ve>0$.
Hence,
$$
\lim_{T\rightarrow\infty}k^*(\ve,T, y_0)\geq k^*(y_0).
$$
Let us prove the opposite inequality. Take $\beta >0 $
 and let $(\gamma', \xi')\in \Omega(y_0) $
be $\beta$-optimal in \eqref{limits-non-ergodic}, that is,
$$
 \int_{Y\times U}k(y,u)\gamma'(dy,du)\leq k^*(y_0) + \beta.
$$
Then
\begin{equation*}
\begin{aligned}
k^*(\ve,T, y_0)&\leq \int_{Y\times U} k(y,u)\gamma'(dy,du)+
\left(\frac{2M}{T}+M\ve\right)\int_{Y\times U} \xi'(dy,du)\\
&\leq  k^*(y_0) + \beta  + \left(\frac{2M}{T}+M\ve\right)\int_{Y\times U}\xi'(dy,du).
\end{aligned}
\end{equation*}
Hence,
$$
\lim_{\ve\downarrow 0,\,T\rightarrow \infty}k^*(\ve, T,y_0)\leq k^*(y_0) .
$$
It follows now from \eqref{GH3} and \eqref{BB1-T-2} that
$$
\limsup_{T\rightarrow\infty}V_T^{\d}(y_0)\le  k^*(y_0),
$$
which implies via \eqref{MC2} that
$$
\limsup_{T\rightarrow\infty} V_T(y_0)\leq k^*(y_0) \ \ \forall \ y_0\in Y.
$$
The part (a) of the theorem is established.  \hf

\begin{remark}
{\rm As mentioned above, the proof of   part {\rm (b)} is similar to that of   part (a). The difference  is that, instead of \eqref{e-before-lim-1}, one needs to use the inequality
\begin{equation}\label{e-before-lim-2}
\int_{Y\times U} h_{\l}(y)\,\g(dy,du)\le   \int_{Y\times U} k(y,u)\,\g(dy,du)\ \ \ \forall \ \gamma\in W.
\end{equation}
Under the assumption (A7), the validity of (\ref{e-before-lim-2}) for any $\l>0$ was established  within the proof of Theorem 3.4 in  \cite{BQR-2015}.}
\end{remark}

\bigskip

From \eqref{limits-non-ergodic-dual-4}, Theorem \ref{Th-upper-bound}, and Proposition \ref{Prop-former-2-3} it follows that under the corresponding assumptions, for any $y_0\in Y$
\begin{equation}\label{BX1}
d^*(y_0)\le \liminf_{T\to \infty} V_T(y_0)\le \limsup_{T\to \infty} V_T(y_0)\le k^*(y_0)
\end{equation}
and
\begin{equation}\label{BX9}
d^*(y_0)\le \liminf_{\l\dn 0} h_{\l}(y_0)\le \limsup_{\l\dn 0} h_{\l}(y_0)\le k^*(y_0).
\end{equation}
These formulas imply that under the strong duality assumption  $d^*(y_0)=k^*(y_0)$ the limits
\begin{equation}\label{lim-exists-1}
V(y_0):=\lim_{T\rightarrow\infty} V_{T}(y_0)
\end{equation}
and
\begin{equation}\label{lim-exists-2}
h(y_0):=\lim_{\l\downarrow 0} h_{\l}(y_0)
\end{equation}
exist and all the inequalities in \eqref{BX1}-\eqref{BX9} hold as equalities. However, strong duality may not hold, and next we investigate situations when individual inequalities in \eqref{BX1}-\eqref{BX9} hold as equalities without  the strong duality assumption. First,  consider the case when there exists an optimal process that is periodic for $t\ge \bar t$ with some $\bar t\ge 0$.

For $\mathcal{T}>0 $ let $(y_{\mathcal{T}}(\cdot), u_{\mathcal{T}}(\cdot)) $ be  a $\mathcal{T} $-periodic admissible process. This process will be referred to as {\it finite time (FT) reachable from $y_0$}  if there exist  $\bar t \geq 0 $ and a control $u(\cdot)\in \U_{\bar t}(y_0)$ such that the  solution $y(t)= y(t,y_0,u) $ of \eqref{e-CSO}
obtained with this control satisfies the equality $y(\bar t) = y_{\mathcal{T}}(0) $.

 Consider the optimal control problem
\begin{equation*}\label{e-main-8-1}
\inf_{\mathcal{T},\left(y_{\mathcal{T}}(\cdot),u_{\mathcal{T}}(\cdot)\right)  }\left\{\frac{1}{\mathcal{T}}\int_{0}^{\mathcal{T}} k(y_{\mathcal{T}}(t),u_{\mathcal{T}}(t))\,dt \right\} :=V_{per}(y_0),
\end{equation*}
where ${\rm inf}$ is over all $\mathcal{T}>0$ and over all $\mathcal{T}$-periodic pairs $(y_{\mathcal{T}}(\cdot), u_{\mathcal{T}}(\cdot)) $ that are FT reachable from $y_0$. Clearly,
\begin{equation}\label{e:occup-meas-def-eq-per-1}
V_{per}(y_0)\geq \liminf_{T\rightarrow\infty}V_T(y_0).
\end{equation}

\bigskip
The following theorem is proved in \cite{BG}, Corollary 3.1.

\begin{Theorem}\label{Prop-former-2-4}
If \eqref{e:occup-meas-def-eq-per-1} holds as equality, that is, $\liminf$ in \eqref{e:occup-meas-def-eq-per-1} is reached on a sequence of periodic processes FT-reachable from $y_0$, then
$$
\liminf_{T\rightarrow\infty}V_T(y_0)\geq k^*(y_0).
$$
\end{Theorem}

This theorem, along with Theorem \ref{Th-upper-bound} (a), implies the following corollary.
\begin{Corollary}\label{Prop-former-2-4-Cor}
Let the assumptions of  Theorem \ref{Th-upper-bound} (a) hold, and $\liminf$ in \eqref{e:occup-meas-def-eq-per-1} is reached on a sequence of periodic processes FT-reachable from $y_0$. Then the limit $\disp V(y_0)=\lim_{T\to \infty}V_T(y_0)$ exists and
$$
V(y_0)=k^*(y_0).
$$
\end{Corollary}

\bigskip

So far we have not assumed that
the limit optimal value functions $V(\cdot)$ and $h(\cdot)$ are differentiable (this is a strong assumption, but it may hold in nontrivial examples, as demonstrated in Section \ref{Examples}). However, if we do assume it, then, as shown in the following theorem, $V(\cdot)$ and $h(\cdot)$
are equal to the optimal value $d^*(y_0)$ of the dual problem \eqref{limits-non-ergodic-dual} without the strong duality assumption $d^*(y_0)=k^*(y_0)$.

\begin{Theorem}\label{ThN1}
{\rm (a)} Let the pointwise limit \eqref{lim-exists-1} exist and $V\in C^1(\bar Y)$, that is, $V(\cdot)$ is continuously differentiable on the interior of $Y$, and $\nabla V$ can be extended by continuity to the boundary of $Y$.  Then
\begin{equation*}\label{CC10}
 V(y_0)= d^*(y_0) \ \ \forall \ y_0\in Y.
\end{equation*}
{\rm (b)} Let the pointwise limit \eqref{lim-exists-2}
exist and  $h\in C^1(\bar Y)$. Then
\begin{equation*}\label{abel-lim-4-1}
 h(y_0)= d^*(y_0) \ \ \ \forall \ y_0\in Y.
\end{equation*}
{\rm (c)} Consequently, if the assumptions of both parts {\rm (a)} and {\rm (b)} hold, then Abel and Ces\`aro limits are equal, that is,
$$
V(y_0)=h(y_0)  \ \ \forall \ y_0\in Y.
$$
\end{Theorem}

The assumption of continuous differentiability of $V$ is essential for the validity of this theorem. This is in contrast to discrete-time systems where the counterpart of Theorem \ref{ThN1} holds
for merely continuous $V(\cdot)$ and $h(\cdot)$  (\cite{BGS}, Theorem 4.2). At the same time, a known sufficient condition ensuring the validity of the equality $V(y_0)=h(y_0)$ for all $y_0\in Y$ is the {\em uniform} convergence of $V_T(\cdot)$ and $h_{\l}(\cdot)$ in  \eqref{lim-exists-1} and \eqref{lim-exists-2} (see Theorems 3.4 and 3.8 in \cite{BQR-2015}). Let us emphasize  that  in
Theorem \ref{ThN1}, the convergence to the limits is not assumed to be uniform.

{\bf Proof of Theorem \ref{ThN1}.}
From \eqref{e-main-1} it follows that $V(y_0)\ge d^*(y_0)$ for all $y_0\in Y$; therefore, to prove part (a), it remains to show the opposite inequality.

From \eqref{e-before-lim-1} it follows that
\begin{equation}\label{e-after-lim}
\min_{\g\in W}\int_{Y\times U} (k(y,u)-V(y))\,\g(dy,du)\ge 0.
\end{equation}
The problem on the left hand side of (\ref{e-after-lim}), i.e.,
\begin{equation}\label{CC8-2}
\min_{\g\in W}\int_{Y\times U}(k(y,u)-V(y))\,\gm,
\end{equation}
is an IDLP problem, whose dual is
\begin{equation}\label{BB12-2}
\sup_{\eta\in C^1}\min_{(y,u)\in Y\times U}\{k(y,u)-V(y)+\nabla \eta(y)f(y,u)\}
\end{equation}
(cf. \eqref{CC8-1}-\eqref{BB12-1}).
Since the optimal values of \eqref{CC8-2} and \eqref{BB12-2} are equal, \eqref{e-after-lim} is equivalent to
\begin{equation*}\label{BB15-2}
\sup_{\eta\in C^1}\min_{(y,u)\in Y\times U}\{k(y,u)-V(y)+\nabla \eta(y)f(y,u)\}\ge 0.
\end{equation*}
Therefore, for any $\ve>0$ there exists a function $\eta_{\ve}(\cdot)\in C^1$ such that
\begin{equation}\label{e-eps-feasib}
k(y,u)-V(y)+\nabla\eta_{\ve}(y)f(y,u)\ge -\ve\quad \hbox{for all }(y,u)\in Y\times U.
\end{equation}
For $y_0\in\,$ int $Y$, similarly to \eqref{XX3} we obtain
\begin{equation*}
TV_T(y_0)\le k(y_0,u)\D t+TV_T(y(\D t))+M\D t+o(\D t)
\end{equation*}
for sufficiently small $\D t$, where $M$ is such that $|k(y,u)|\le M$ for all $(y,u)\in Y$.
Dividing both sides of this inequality by $T$ and then passing to the limit as $T\to \infty$, we obtain
\begin{equation}\label{XV1}
V(y_0)\le V(y(\D t))
\end{equation}
Under the differentiability assumpton, the latter implies that
\begin{equation}\label{PZ2}
\nabla V(y_0)f(y_0,u)\ge 0 \hbox{ for all }(y_0,u)\in ({\rm int}\, Y)\times U.
\end{equation}
Due to the assumption on continuity of $\nabla V$ on $Y$, inequality above holds for all $y_0\in Y$.

Consider the problem
\begin{equation}\label{BB21-2}
\sup_{(\psi,\eta)\in Q} \psi(y_0),
\end{equation}
where $Q$ is the set of pairs $(\psi(\cdot),\eta(\cdot))\in C^1\times C^1$ that satisfy the inequalities
\begin{equation}\label{BB25-2}
\begin{aligned}
&k(y,u)-\psi(y)+\nabla \eta(y)f(y,u)\ge 0,\\
&\nabla \psi(y)f(y,u)\ge 0\ \ \ \ \forall \ (y,u)\in Y\times U.
\end{aligned}
\end{equation}
We will show below that the optimal value in \eqref{BB21-2}-\eqref{BB25-2} is equal to $d^*(y_0)$. Assuming that this is true, due to Whitney Extension Theorem (\cite{Whitney}), function $V$ can be extended from $Y$ to $\reals^n$ as a $C^1$ function. Since the pair $\psi=V-\ve$ and $\eta=\eta_{\ve}$ satisfies \eqref{BB25-2} due to \eqref{e-eps-feasib} and \eqref{PZ2}, we conclude that
\begin{equation*}\label{ZZ10a}
V(y_0)-\ve\le d^*(y_0),
\end{equation*}
which implies the validity of part (a) of Theorem \ref{ThN1}.

Let us now prove the above claim that  the optimal value in \eqref{BB21-2}-\eqref{BB25-2} is equal to $d^*(y_0)$.
Let us denote the value of supremum in \eqref{BB21-2} by $\hd^*(y_0)$. The inequality $\hd^*(y_0) \leq  d^*(y_0) $ is true, since for any pair $\ (\psi(\cdot),\eta(\cdot))\in Q$, the triplet $\ (\mu , \psi(\cdot),\eta(\cdot))\in \Dd $ with $\mu = \psi(y_0) $ and where $\Dd$ is given after formula \eqref{limits-non-ergodic-dual}. Let us prove the opposite inequality. Let a triplet $\ (\mu' , \psi'(\cdot),\eta'(\cdot))\in \Dd$ be such that $\mu'\geq d^*(y_0)-\delta $, with $\delta > 0 $ being arbitrarily small. Then the pair $\ (\tilde \psi'(\cdot), \eta'(\cdot))\in Q$, with $\tilde\psi'(y)= \psi'(y)- \psi'(y_0) + \mu'$. Since $ \tilde\psi'(y_0) = \mu'$, it leads to the inequality $\hat d^*(y_0) \geq  d^*(y_0) - \delta $ and, consequently, to the inequality $\hat d^*(y_0) \geq  d^*(y_0) $ since $\delta > 0 $ is arbitrarily small. Thus, the claim and part (a) of the theorem are proved.  The proof of part (b) is similar and is omitted.\hf

\section{Sufficient Optimality Conditions}\label{opt_cond}

In this section we consider the optimal control problem
\begin{equation}\label{CC21}
\inf_{u(\cdot)\in \U(y_0)} \liminf_{T\to \infty} {1\o T}\int_0^T k(y(t),u(t))\,dt
\end{equation}
and its sufficient optimality conditions in terms of the maximizing functions in the dual problem. First we relate the value function of (\ref{CC21})  with that of \eqref{Cesaro}.

\begin{Proposition}\label{PN3} The optimal value in problem  (\ref{CC21}) is equal to $\disp\liminf_{T\to \infty} V_T(y_0) $. That is,
for all $y_0\in Y$
$$
\inf_{u(\cdot)\in \U(y_0)} \liminf_{T\to \infty} {1\o T}\int_0^T k(u(t),y(t))\,dt=\liminf_{T\to \infty} V_T(y_0).
$$
\end{Proposition}
\bigskip

{\bf Proof.}
Let $u(\cdot)\in \U(y_0)$ and let $y(\cdot)=y(\cdot,y_0,u(\cdot))$ be the corresponding trajectory. Then
$$
{1\o T}\int_0^T k(y(t),u(t))\,dt\ge V_T(y_0).
$$
Therefore,
$$
\liminf_{T\to \infty}{1\o T}\int_0^T k(y(t),u(t))\,dt\ge \liminf_{T\to \infty} V_T(y_0)
$$
and, hence,
$$
\inf_{u(\cdot)\in \U(y_0)} \liminf_{T\to \infty} {1\o T}\int_0^T k(y(t),u(t))\,dt\ge\liminf_{T\to \infty} V_T(y_0).
$$
Let us prove the opposite inequality.
For any $\varepsilon>0$, $u(\cdot)\in \U(y_0)$,  and for sufficiently large $T$,
$$
{1\o T}\int_0^T k(u(t),y(t))\,dt \ge \liminf_{T\to \infty} {1\o T}\int_0^T k(u(t),y(t))\,dt -\varepsilon ,
$$
where $y(\cdot)=y(t,y_0,u) $.
Therefore,
\begin{equation*}
\begin{aligned}
{1\o T}\int_0^T & k(u(t),y(t))\,dt \ge
\inf_{u'\in \U(y_0)}\liminf_{T\to \infty} {1\o T}\int_0^T k(u'(t),y'(t))\,dt -\varepsilon
\end{aligned}
\end{equation*}
and, consequently,
\begin{equation*}
\begin{aligned}
V_T(y_0) \ge \inf_{u\in \U(y_0)}\liminf_{T\to \infty} {1\o T}\int_0^T k(u(t),y(t))\,dt -\varepsilon.
\end{aligned}
\end{equation*}
Hence,
$$
\liminf_{T\to \infty}V_T(y_0)\ge \inf_{u\in \U(y_0)}\liminf_{T\to \infty} {1\o T}\int_0^T k(u(t),y(t))\,dt.
$$
The proposition is proved. \hf

\bigskip

\begin{Theorem}\label{Th-suff}
Assume that a pair $(\bp,\be)\in C^1\times C^1$ delivering maximum in problem \eqref{limits-non-ergodic-dual}-\eqref{limits-non-ergodic-dual-2} (or \eqref{CC9}) exists and, for some  admissible process $(y(\cdot),u(\cdot))$ and all $t\ge 0$,
\begin{equation}\label{XX6}
\begin{aligned}
k(y(t),u(t))+(\bp(y_0)-\bp(y(t)))+
\nabla\bar\eta(y(t))f(y(t),u(t))=d^*(y_0).
\end{aligned}
\end{equation}
Then\\
{\rm (a)} there exists the limit $\disp V(y_0)=\lim_{T\rightarrow\infty}V_T(y_0)$;\\
{\rm (b)} there holds the equality
\begin{equation}\label{XX7}
\disp V(y_0)= d^*(y_0);
\end{equation}\\
{\rm (c)} the process $(y(\cdot),u(\cdot))$ is optimal in (\ref{CC21}).
\end{Theorem}

{\bf Proof.} Integrating \eqref{XX6}  and taking into account that $\bp(y_0)-\bp(y(t))\le 0$ due to
\eqref{limits-non-ergodic-dual-2}, we obtain
$$
{1\o T}\int_0^T k(y(t),u(t))\,dt+{1\o T}(\be(y(T))-\be(y_0))\le d^*(y_0).
$$
Since $\disp V_T(y_0)\le {1\o T}\int_0^T k(y(t),u(t))\,dt$ and the second term vanishes as $T\to \infty$, the latter implies that
$$
\limsup_{T\to \infty} V_T(y_0)\le d^*(y_0).
$$
Taking into account \eqref{e-main-1}, we conclude that the limit $\disp V(y_0)=\lim_{T\to \infty} V_T(y_0)$ exists
and is equal to $d^*(y_0)$, that is, parts (a) and (b) of the theorem are true. We also obtain that\\
$\disp \lim_{T\to \infty}{1\o T}\int_0^T k(y(t),u(t))\,dt=d^*(y_0)$, and due to \eqref{XX7}, we conclude that
$$
\lim_{T\to \infty}{1\o T}\int_0^T k(y(t),u(t))\,dt=V(y_0).
$$
Therefore, the process $(y(\cdot),u(\cdot))$ is optimal, and part (c) is true. The theorem is proved. \hf

\bigskip

An important corollary of Theorem \ref{Th-suff} is the following proposition.

\begin{proposition}\label{Rem-feedback}
Under the assumptions of Theorem \ref{Th-suff}
\begin{equation*}
\begin{aligned}
(y(t),u(t))=
{\rm argmin}_{(y,u)\in Y\times U}\{k(y,u)-\bar\psi(y)+\nabla\bar\eta(y)f(y,u)\} \quad\forall\,t\ge 0,
\end{aligned}
\end{equation*}
therefore,
\begin{equation}\label{feedback}
u(t)={\rm argmin}_{u\in U}\{k(y(t),u)+\nabla\bar\eta(y(t))f(y(t),u)\}\,.
\end{equation}
The latter implies the optimal feedback control law
\begin{equation}\label{e-feedback}
\begin{aligned}
u^f[y]={\rm argmin}_{u\in U}\{k(y ,u)+\nabla\bar\eta(y)f(y,u)\}.
\end{aligned}
\end{equation}
\end{proposition}
{\bf Proof.} This proposition follows immediately from \eqref{XX6} and \eqref{CC9}.

\bigskip

The law \eqref{e-feedback} can be used to construct optimal control numerically, which can be a subject of further research.
An example illustrating Proposition \ref{Rem-feedback} is presented in Section \ref{Examples}.

In the situation when the function $V(\cdot)$ is continuously differentiable, in addition to Theorem \ref{Th-suff}, the following is true.

\begin{Theorem}\label{PN2}
Assume
that the pointwise limit $\disp V(\cdot)=\lim_{T\rightarrow\infty}V_T(\cdot)$ exists, is continuously differentiable on the interior of $Y$, and $\nabla V$ can be extended by continuity to the boundary of $Y$.  Then

{\rm (a)}  functions $(\bar\psi,\bar\eta)\in C^1\times{C^1}$ are maximizers in  \eqref{CC9} if and only if $\bar\psi$ satisfies \eqref{limits-non-ergodic-dual-2} and
 \begin{equation}\label{min}
\begin{aligned}
\min_{(y,u)}\{k(y,u)-\bp(y)+\nabla\bar\eta(y)f(y,u)\}=V(y_0)-\bp(y_0);
\end{aligned}
\end{equation}

{\rm (b)} if $\bar\eta\in {C^1}$ is such that
\begin{equation}\label{CC7}
\min_{(y,u)}\{k(y,u)-V(y)+\nabla\bar\eta(y)f(y,u)\}= 0,
\end{equation}
then the functions $\psi=V$ and $\eta=\bar\eta$ are maximizers in \eqref{CC9}.


\end{Theorem}

{\bf Proof.} For  $(\bar\psi,\bar\eta)\in C^1\times{C^1}$ being maximizers in  \eqref{CC9} is equivalent to
 \begin{equation}\label{min2}
\min_{(y,u)}\{k(y,u)+(\bp(y_0)-\bp(y))+\nabla\bar\eta(y)f(y,u)\}=d^*(y_0).
\end{equation}
Due to Theorem \ref{ThN1} (a) we have $d^*(y_0)=V(y_0)$, which implies part (a). To show the validity of (b) notice that if $\bar\eta$ is such that \eqref{CC7} holds, then the  functions $\psi=V$ and $\eta=\bar\eta$ satisfy \eqref{min}.  Due to part (a) of the theorem, these functions are maximizers in \eqref{CC9}.
Theorem is proved. \hf

\begin{Remark}\label{R-non-unique}
{\rm Due to Theorem \ref{PN2}{ (b)}, $\psi=V$, along with appropriate $\eta$, is a maximizer in \eqref{CC9}. However, the set of maximizing $\psi$ may be much larger then $V$ alone, as is demonstrated by the example in the following section.}
\end{Remark}

\section{An Example}\label{Examples}

In this section we illustrate an application of  Proposition \ref{Rem-feedback} and Theorem \ref{PN2} with an an example.
Consider the system (Example 5.1 in \cite{BG})
\begin{equation}\label{VC1}
\begin{aligned}
&y_1'=y_2u,\\
&y_2'=-y_1u,\\
&(y_1(0),y_2(0))=(y_{10},y_{20}),\\
&u\in [-1,1]\,
\end{aligned}
\end{equation}
with the state constraint
\begin{equation}\label{VM4}
y_{1}^2+y_{2}^2\le 1\,.
\end{equation}
Since $\disp {d\o dt}(y_1^2+y_2^2)=0$ and $(y_1')^2+(y_2')^2=(y_{01}^2+y_{02}^2)u^2$, this system describes rotation about the origin with the angular speed $u$.
The goal is to find $u$ that minimizes the long-run average distance from the point $(1,0)$:
\begin{equation}\label{VX3}
V(y_{01},y_{02}):=\min_{u(\cdot)}\liminf_{T\to \infty} {1\o T}\int_0^T ((1-y_1)^2+y_2^2)\,dt,
\end{equation}
that is, $k(y_1,y_2,u)=(1-y_1)^2+y_2^2$.

For convenience of analysis, write  \eqref{VC1} in polar coordinates:
\begin{equation}\label{e-extra}
\begin{aligned}
&r'=0,\\
&\t'=u,\\
& (r(0),\t(0))=(r_0,\t_0),\\
&u\in [-1,1]\,.
\end{aligned}
\end{equation}
The constraint \eqref{VM4} becomes
\begin{equation}\label{VM5}
0\le r\le 1,\;-\pi\le \t\le \pi.
\end{equation}
Note that the constraint (\ref{VM4}) is equivalent to just $0\le r\le 1 $. The introduction of the additional inequality $-\pi\le \t\le \pi$ in (\ref{VM5}) makes the state constraint set compact in polar coordinates without changing
the optimal value of the problem. Note also 
 that system \eqref{VC1} is invariant on the set \eqref{VM4}, and system (\ref{e-extra}) is viable (but not invariant) on the set \eqref{VM5}.

We have $(1-y_1)^2+y_2^2=1-2r\cos\t+r^2$,
and we denote
\begin{equation}\label{BB2}
V_T(r_0,\t_0):={\min_{u(\cdot)}}\,{1\o T}\int_0^T (1-2r\cos\t+r^2)\,dt.
\end{equation}
Due to Proposition \ref{PN3}, the value function in \eqref{VX3} is equal to $\lim_{T\to \infty} V_T.$

It is clear that if $\t_0\in(0,\pi]$, the optimal control is $u=-1$ until $\t=0$, then $u=0$;  if $\t_0\in[-\pi,0)$, the optimal control is $u=1$ until $\t=0$, then $u=0$.
The time when $\t$ reaches 0 is equal to $|\t_0|$. Therefore, if $\t_0\in(0,\pi]$, then, along the optimal process,
$$
\t(t)=\begin{cases}
\t_0-t,\,& t\in[0,\t_0]\\
0,&t>\t_0
\end{cases}
$$
and, for $T>\t_0$,
\begin{equation}\label{BX3}
V_T(r_0,\t_0)={1\o T}\left(\int_0^{\t_0}(1-2r\cos\t+r^2)\,dt+(T-\t_0)(1-r_0)^2\right).
\end{equation}
If $\t_0\in[-\pi,0)$, then
$$
\t(t)=\begin{cases}
\t_0+t,\,& t\in [0,-\t_0]\\
0,&t>-\t_0\,
\end{cases}
$$
and, for $T>-\t_0$,
\begin{equation}\label{BB5}
V_T(r_0,\t_0)={1\o T}\left(\int_0^{-\t_0}(1-2r\cos\t+r^2)\,dt+(T+\t_0)(1-r_0)^2\right).
\end{equation}
For $\t_0\in (0,\pi]$ we have
$$
\int_0^{\t_0}\cos\t(t)\,dt=\int_0^{\t_0}\cos(\t_0-t)\,dt=\sin \t_0.
$$
Taking into account that $r\equiv r_0$, from \eqref{BX3} we get
\begin{equation}\label{BB4}
V_T(r_0,\t_0)={1\o T}(\t_0(1+r_0^2)-2r_0\sin \t_0+(T-\t_0)(1-r_0)^2)=(1-r_0)^2+{2r_0\o T}(\t_0-\sin\t_0)\,.
\end{equation}
For $\t_0\in [-\pi,0)$ we have
$$
\int_0^{-\t_0}\cos\t\,dt=\int_0^{-\t_0}\cos(\t_0+t)\,dt=-\sin \t_0,
$$
and, from \eqref{BB5},
\begin{equation}\label{BB6}
V_T(r_0,\t_0)=(1-r_0)^2+{2r_0\o T}(-\t_0+\sin\t_0).
\end{equation}
Formulas \eqref{BB4} and \eqref{BB6} can be combined as
$$
V_T(r,\t)=(1-r)^2+{2r\o T}|\t-\sin\t|,\;\t\in[-\pi,\pi].
$$
From this formula we obtain
\begin{equation*}
V(r,\t)=\lim_{T\to \infty}V_T(r,\t)=(1-r)^2.
\end{equation*}
This function is of class $C^1$.
Due to Theorem \ref{PN2} (b), if $\eta\in C^1$ is such that
\begin{equation}\label{VC2}
\min_{r,\t,u} \{k(r,\t)-V(r,\t)+\nabla\eta(r,\t)f(r,\t,u)\}=0,
\end{equation}
then $\psi=V$ and $\eta$ are maximizing functions in the dual problem. Let us verify that
\begin{equation}\label{e-eta}
\eta(r,\t)=2r|\t-\sin\t|
\end{equation}
satisfies \eqref{VC2}. First, notice that $\eta$ is continuously differentiable, since both one-sided derivatives at $\t=0$ are equal to zero. (In fact, $\eta\in C^2$, but the third derivative is discontinuous at $\t=0$.) We have
\begin{equation}\label{BB7}
\begin{aligned}
&\nabla \eta(r,\t)f(r,\t,u)=2r\,{\rm sgn}\,\t \,(1-\cos\t)u,\\
& \min_{u\in [-1,1]} \{\nabla\eta(r,\t)f(r,\t,u)\}=-2r(1-\cos\t)
\end{aligned}
\end{equation}
and
$$
k(r,\t)-V(r,\t)+\min_{u\in [-1,1]} \nabla\eta(r,\t)f(r,\t,u)= (1-2r\cos\t+r^2)-(1-r)^2-2r(1-\cos\t)\equiv 0.
$$
Therefore, \eqref{VC2} holds. Consequently,  $\psi(r)=V(r)=(1-r)^2$ and    $\eta(r,\t)$ defined by (\ref{e-eta}) are maximizing functions in the dual problem \eqref{CC9}.


Note that the application of \eqref{e-feedback} in \eqref{BB7} leads to
the feedback law
$$
u^f[\t]=\begin{cases} 1,&\t<0\\ -1,&\t>0\end{cases}
$$
consistent with the observation made after formula \eqref{BB2}.

As mentioned in Remark \ref{R-non-unique}, $\psi=V$ (along with $\eta$ given by \eqref{e-eta}) is not the only maximizer in \eqref{CC9}. Indeed,
due to Theorem \ref{PN2}(a), functions $(\psi, \eta)$ are maximizers in \eqref{CC9} if and only if $\psi$ satisfies
\begin{equation}\label{VM2}
\nabla\psi(r,\t)f(r,\t,u)\ge 0 \quad \hbox{for all }r,\t,u
\end{equation}
and
\begin{equation}\label{VM3}
\min_{(r,\t,u)}\{k(r,\t)-\psi(r,\t)+\nabla\eta(r,\t)f(r,\t,u)\}=V(r_0,\t_0)-\psi(r_0,\t_0).
\end{equation}
Since $f(r,\t,u)=(0,u)$, relation \eqref{VM2} becomes  $\disp {\pl \psi\o \pl \t}u\ge 0$  for all $u\in [-1,1]$. Hence, $\disp {\pl \psi\o \pl \t}=0$, that is, $\psi$ may not depend on $\t$. Since
$$
k(r,\t)-\psi(r)+\min_{u\in [-1,1]} \nabla\eta(r,\t)f(r,\t,u)= (1-2r\cos\t+r^2)-\psi(r)-2r(1-\cos\t) =(1-r)^2-\psi(r),
$$
relation \eqref{VM3} becomes
$$
\min_{r\in [0,1]}\{(1-r)^2-\psi(r)\}=(1-r_0)^2-\psi(r_0).
$$
Thus, any $\psi(r)$ such $\disp\min_{r\in [0,1]} \{(1-r)^2-\psi(r)\}$ is reached at $r_0$, along with $\eta$ given by \eqref{e-eta}, is a maximizer in \eqref{CC9}.
A specific example is $r_0=1$ and $\psi(r)=(1-r)^3$.

Email addresses of the authors: 

V. Gaitsgory vladimir.gaitsgory@mq.edu.au\\
I. Shvartsman ius13@psu.edu (corresponding author)

\end{document}